\documentclass[a4paper,12pt,reqno,twosides,onecolumn,final]{amsart}

\usepackage{amsmath,amsfonts,amssymb,amscd}

\usepackage[foot]{amsaddr} 

\usepackage{graphicx}
\usepackage{epstopdf}

\usepackage{bm}  
\usepackage[left,pagewise]{lineno}

\usepackage{setspace}

\usepackage[utf8]{inputenc} 
\usepackage{lmodern}
\usepackage{subfig}
\usepackage{natbib}

\usepackage[]{algorithm}
\usepackage[]{algorithmic}
\usepackage{amssymb}
\usepackage{mathabx}

\usepackage{color}
\usepackage{caption}

\usepackage[skip=5pt,indent=15pt]{parskip}

\usepackage{chngcntr}

\newcommand{\av}[1]{\ensuremath{\left\{\!\!\left\{#1\right\}\!\!\right\}}}
\newcommand{\jp}[1]{\ensuremath{\left[\!\left[#1\right]\!\right]}}

\DeclareMathOperator{\Div}{\mathrm{div}}
\DeclareMathOperator{\Kappa}{\mathcal{K}}

\let\div\relax

\DeclareMathOperator{\div}{div}
\DeclareMathOperator{\curl}{\bf{curl}}

\def\0{\phantom{0}}

\newtheorem{theorem}{Theorem}[section]

\theoremstyle{definition}

\theoremstyle{remark}

\begin{document}


\title[On FEM for Heterogeneous Elliptic Problems]{On Finite Element Methods for Heterogeneous Elliptic Problems$^\star$
}



\author[A.F.D. Loula]{Abimael F. D. Loula}
\address[A.F.D. Loula, J.N.C. Guerreiro, E.M. Toledo]{LNCC - Laboratório Nacional de Computação Científica, Petrópolis, RJ, Brasil}
\curraddr{}
\email[A.F.D. Loula]{aloc@lncc.br}
\thanks{}

\author[M.R. Correa]{Maicon R. Correa}
\address[M.R. Correa]{Universidade Estadual de Campinas (UNICAMP), Departamento de Matemática Aplicada, IMECC, Brasil}
\email[M.R. Correa]{maicon@ime.unicamp.br}
\thanks{}

\author[J.N.C. Guerreiro]{João N. C. Guerreiro}
\curraddr{}
\email[J.N.C. Guerreiro]{joao@lncc.br}
\thanks{}

\author[E.M. Toledo]{Elson M. Toledo}
\curraddr{}
\email[E.M. Toledo]{emtc@lncc.br}
\thanks{}


\date{}

\dedicatory{$^\star$Dedicated to Luiz Bevilaccqua, Raul Feij\'oo and Luis Rojas.
\\ \vspace{.5cm}
{\rm This is a preprint of the original article published in} \\{\it International Journal of Solids and Structures, 45:1 (2008), pp 6436-6450}. \\
\rm https://doi.org/10.1016/j.ijsolstr.2008.08.005
}

\begin{abstract}
Dealing with variational formulations of second order elliptic
problems with discontinuous coefficients, we recall a single field
minimization problem of an extended functional presented by
\cite{BEVI74}, which we associate with the basic idea supporting
discontinuous Galerkin finite element methods. We review residual
based stabilized mixed methods  applied to Darcy flow in homogeneous
porous media and extend them to heterogeneous media with an
interface of discontinuity. For smooth interfaces, the proposed
formulations preserve the continuity of the flux and exactly imposes
the constraint between the tangent components of Darcy velocity on
the interface.
Convergence studies for a heterogeneous and anisotropic porous
medium confirm the same rates of convergence predicted for
homogeneous problem with smooth solutions.
\end{abstract}


\maketitle



\section{Introduction}

We dedicate this paper to Luiz Bevilacqua, Raul Feij\'oo, and Luis
Rojas after their work {\it ``A variational principle for the
Laplace operator with application in the torsion of composite rods"}
\citep{BEVI74}, published in International Journal of Solids and
Structures in its ninth year of existence after being created by
George Herrmann in 1965. It is also worth mentioning that at that
time Brazil had graduated just one PhD in all branches of
engineering, at COPPE/UFRJ (now named Instituto Alberto Luiz Coimbra
in honor to its founder) where that work was developed. In that
paper the authors present a Ritz method with discontinuous
coordinate functions to minimize an extended functional with the
Erdman-Weierstrass corner conditions \citep{GELFAND63} weakly
imposed.

Variational formulations in broken function spaces with continuity
naturally imposed in a weak sense, as in \cite{BEVI74}, have become
very popular in the context of discontinuous Galerkin finite element
methods \citep{DDUPONT76,WHEELER78,ARNOLD82,ODEN98,DUARTE2000,
DUTRA02,BREZZI2005,DAWSON2005,SULI05,BREZZI2006,ALVAREZ2006,HUGHES2006A}.
We acknowledge this fact in dealing with Darcy's flow in
heterogeneous porous media with an interface of discontinuity. This
problem basically consists of the mass conservation equation plus
Darcy's law, that relates the average velocity of the fluid in a
porous medium with the gradient of a potential field through the
hydraulic conductivity tensor.

The simplest way to solve Darcy's system is based on the second
order elliptic problem obtained by substituting Darcy's law in the
mass conservation equation leading to a Poisson equation in the
potential field. After solving this Poisson problem for the
potential, velocity is calculated by taking the gradient of the
solution multiplied by the hydraulic conductivity. Constructing
finite element approximations based on this formulation is
straightforward. However, this direct approach leads to lower-order
approximations for velocity compared to potential and, additionally,
the corresponding balance equation is enforced in an extremely weak
sense. Alternative formulations have been developed to enhance the
velocity approximation, like post-processing techniques
\citep{CORDES92,LOULA95,RITA2002,CORREA2007} and mixed methods
\citep{RAVIART77,BREZZI85,MOSE94,DURLOFSKY94}. The basic idea of
post-processing formulations is to use the optimal stability and
convergence properties of the classical Galerkin approximation for
the potential field to derive more accurate velocity approximations
than that obtained by direct application of Darcy's law.

Mixed methods are based on the simultaneous approximation of
potential and velocity fields, and their main characteristic is the
use of different spaces for velocity and potential to satisfy a
compatibility condition between the finite element spaces, which
reduces the flexibility in constructing finite approximations
\citep[LBB condition - ][]{BREZZI74,FORTIN91}. A well known
successful approach is the dual mixed formulation de\-ve\-lo\-ped by
\cite{RAVIART77} using divergence based finite element spaces for
the velocity field combined with discontinuous Lagrangian spaces for
the potential. Stabilized mixed finite element methods have been
proposed to overcome the compatibility conditions between the finite
element spaces, see
\cite{LOULA87A,FRANCA88,FRANCA88A,LOULA88,ELSON90,MASUD2002,BREZZI2005,BARRENECHEA2007,CORREA2008A}
and references therein.
In general, stabilized formulations use continuous Lagrangian finite
element spaces and have been successfully employed in simulating
Darcy flows in homogeneous porous media. The point is their
applicability to problems where the porous formation is composed of
subdomains with different conductivities. On the interface between
these subdomains, the normal component of Darcy velocity must be
continuous (mass conservation) but the tangential component is
discontinuous, and {any} formulation based on $C^0(\Omega)$
Lagrangian interpolation for velocity fails  in representing the
tangential discontinuity, producing inaccurate approximations and
spurious oscillations, while completely discontinuous approximations
for velocity do not necessarily fit the continuity of the flux.
Possible alternatives are the discontinuous Galerkin methods
\citep[e.g.][]{BREZZI2005,HUGHES2006A}, or coupling continuous and
discontinuous Galerkin methods \citep[e.g.][]{DAWSON2005}.

In the present work we study finite element approximations for
second order elliptic problems with variable coefficients presenting
an interface of discontinuity. As a model problem we consider Darcy
flow in heterogeneous porous media composed of layers with different
conductivities.
The model problem is established in Section \ref{sec:model}.
Compatible $C^0(\Omega)$ Lagrangian formulations, including the
classical Galerkin method and Stabilized mixed methods are reviewed
 in Section \ref{sec:aprox}, testing their behavior when applied to
 heterogeneous problems.
In Section \ref{sec:broken} we present weak formulations in broken
spaces and relate the extended functional introduced by
\cite{BEVI74} with Discontinuous Galerkin methods.
In Section \ref{sec:discont} we propose stabilized mixed
formulations with selective equal-order continuous/discontinuous
Lagrangian finite element spaces adopted for both velocity and
potential fields. The discontinuity of the tangential component of
the velocity field is appropriately captured and the continuity of
the potential and flux are strongly imposed on the interface of the
two media as proposed in \citet{CORREA2006} and \citet{CORREA2007}.
Finally, in Section \ref{sec:conclu} we draw some conclusions.
%

%
%
\section{Model Problem}
\label{sec:model}
%
%

\begin{figure}[htb]
\centering
\includegraphics[width=.5\linewidth]{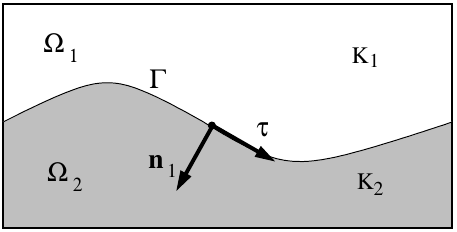}
\caption{Heterogeneous media with an interface of discontinuity.}
\label{fig:meios}%
\end{figure}

Let $\Omega = \Omega_1\cup\Omega_2 \subset \mathbb{R}^2$, with
smooth boundary $\partial \Omega=\Gamma_1 \cup \Gamma_2$ and outward
unit normal vector $\bm{n}$, be the domain of a rigid porous media
composed of two subdomains $\Omega_1$ and $\Omega_2$ separated by a
smooth interface $\Gamma$ and  saturated with an incompressible
homogeneous fluid as illustrated in Figure \ref{fig:meios}.

For simplicity we will consider the porous media homogeneous and
anisotropic on each subdomain. The mass balance on each subdomain is
established as
\begin{equation}
  \label{eq:balanco1}
  \Div \bm{u}_1 =f_1 \quad \mbox{ in } \quad\Omega_1
\end{equation}
\begin{equation}
  \label{eq:balanco2}
  \Div \bm{u}_2 =f_2 \quad \mbox{ in } \quad\Omega_2
\end{equation}
where $f_i = f|_{\Omega_i}$ are sources of fluid and
$\bm{u}_i=\bm{u}|_{\Omega_i}$, $i=1,2$ are average velocities of the fluid
in each domain $\Omega_i$, given by Darcy's law
\begin{equation}
  \label{eq:darcy1}
  \bm{u}_1=-\Kappa_1\nabla p_1 \quad \mbox{ in } \quad \Omega_1
\end{equation}
\begin{equation}
  \label{eq:darcy2}
  \bm{u}_2=-\Kappa_2\nabla p_2 \quad \mbox{ in } \quad \Omega_2
\end{equation}
that relates the velocity field with the gradient of a hydraulic
potential $p_i$ through the hydraulic conductivity tensor $\Kappa_i=
\Kappa|_{\Omega_i}$. Since $\Kappa_i$ is a bounded and positive definite
tensor, its inverse $\Lambda_i=\Kappa_i^{-1}$, referred to as hydraulic
resistivity tensor, is well defined and also positive definite. We
limit our presentation to homogeneous Dirichlet and Neumann boundary
conditions
  \begin{equation}
   \label{eq:diric}
   p_i = 0  \quad \mbox{ on }  \quad \Gamma_i^\mathrm{d} = \partial \Omega_i^\mathrm{d}/\Gamma ,
  \end{equation}
\begin{equation}
   \label{eq:neum}
   \Kappa_i \nabla p_i \cdot \bm{n} = 0  \quad \mbox{ on }  \quad \Gamma_i^\mathrm{n} = \partial \Omega_i^\mathrm{n}/\Gamma ,
  \end{equation}
with
\( \Gamma_i^\mathrm{d}\cup \Gamma_i^\mathrm{n}=\Gamma_i=\partial
\Omega_i /\Gamma,
\)
\( \; \Gamma_i^\mathrm{d}\cap \Gamma_i^\mathrm{n} = \emptyset\) and
the indexes $``\mathrm{d}"$ and $``\mathrm{n}"$ indicating the part
of the boundary where Dirichlet or Neumann boundary conditions are
prescribed.
 Our model problem consists in solving the
 system  of partial differential equations (\ref{eq:balanco1})--(\ref{eq:darcy2})
 with homogeneous boundary
 conditions (\ref{eq:diric}) and (\ref{eq:neum}), plus the following interface conditions:
\begin{equation}
  \label{eq:nflx}
      \bm{u}_1\cdot \bm{n} = \bm{u}_2\cdot \bm{n}
      \; {\mbox{ on }} \; \Gamma\,,
\end{equation}
\begin{equation}
  \label{eq:conp}
      p_1=p_2   \; {\mbox{ on }} \; \Gamma \,,
\end{equation}
where we adopt $\bm{n}=\bm{n}_2=-\bm{n}_1$ on the interface $\Gamma$.

%
%
\section{Compatible $C^0$ Lagrangian Formulations}
\label{sec:aprox}
%
%

In this Section we present some compatible formulations for the
model problem which are naturally stable, as the classical single
field Galerkin formulation, or stabilized mixed formulations. To
introduce these formulations, we first present some definitions and
notations. Let $L^2(\Omega)$ be the space of square integrable
functions in $\Omega$ with inner product $(\cdot,\cdot)$ and norm
$\|\cdot\|$. $H^m(\Omega)$, for $m\geq 0$ integer, is the Hilbert
space in $\Omega$ with $m$ derivatives in $L^2(\Omega)$ with inner
product $\left( \cdot,\cdot\right)_m$ and norm
$\left\|\cdot\right\|_{m}$. We also define the Hilbert space
\[
H(\div,\Omega)= \left\{\bm{u} \in \left[ L^2(\Omega) \right]^2; \; \Div \bm{u} \in
L^2(\Omega)\right\}
\]
with norm
\[
\|\bm{u} \|_{H(\div)} \coloneq \left\{\| \bm{u} \|^2 + \|\Div \bm{u} \|^2\right \}^{1/2}
\]
and its subspace
\[
H_0(\div,\Omega) = \left\{ \bm{u} \in H(\div,\Omega) ; \bm{u}\cdot\bm{n}=0 \mbox{ on }
\partial \Omega \right\}.
\]
To simplify our presentation,  we consider only essential boundary
conditions.

%
%
\subsection{Standard Galerkin Formulation}
%
%

The simplest approach to our model problem consists in substituting
Darcy's law (\ref{eq:darcy1}) and (\ref{eq:darcy2}) into the mass
conservation (\ref{eq:balanco1}) and (\ref{eq:balanco2}),
respectively, yielding to 

\vspace{.25cm}
\paragraph{Problem P} {\it
  Given the hydraulic conductivity  $\Kappa=\Kappa_i$ in $\Omega_i$ and a function $f$,
  find the potential field $p$ such that
  \begin{equation}  \label{eq:p}
    -\Div (\Kappa \nabla p)   = f   \quad \mbox{ in } \quad \Omega =
    \Omega_1 \cup \Omega_2
  \end{equation}
  with homogeneous boundary condition
  \begin{equation}
    p  = 0  \quad \mbox{ on }  \quad \partial \Omega = \Gamma_1\cup
    \Gamma_2
  \end{equation}
  plus the interface conditions on $\Gamma$.
}

To set the classical primal variational formulation, let
\begin{equation}
\label{eq:q} \mathcal{P} = H^1_0(\Omega)=\{q\in H^1(\Omega) ; q=0 \;
\mbox{on} \;
\partial \Omega \}.
\end{equation}
Multiplying (\ref{eq:p}) by a function $q\in \mathcal{P}$ and integrating it
by parts over $\Omega$, we have: {\it
  Find $p \in \mathcal{P}$ such that
  \begin{equation}
    a(p,q)= f(q)
    \quad \forall\; q \in \mathcal{P}
  \end{equation}
  with
   \begin{equation}
  a(p,q)=  (\Kappa \nabla p,\nabla q)
  \, ; \quad
  f(q)= (f,q).
   \end{equation}  %
}
Considering $f\in H^{-1}(\Omega)$ we have the continuity of
$f(\cdot)$. As $a(\cdot,\cdot)$ is a continuous and $\mathcal{P}-$elliptic
bilinear form, this problem has a unique solution by Lax-Milgram
Lemma \citep{CIARLET78}.

Let  $\{ {\mathcal T}_h \}$ be a family of partitions ${\mathcal
T}=\{\Omega^e\}$ of $\Omega$ indexed by the parameter $h$
representing the maximum diameter of the elements $\Omega^e \in{\mathcal
T}$.
Let $\mathcal{S}_{h}^{k}$ be the $C^0(\Omega)$ Lagrangian finite element space of
degree $k\geq 1$ in each element $\Omega^e$,
\begin{equation}
\label{eq:shk} \mathcal{S}_{h}^{k} = \{ \varphi_h \in C^0(\Omega); \left.
\varphi_h\right|_{\Omega^e} \in \mathbb{P}_k(\Omega^e) \}
\end{equation}
where ${\mathbb{P}}_k(\Omega^e)$ is the set of the polynomials of
degree $\leq k$ posed on $\Omega^e$.
Defining the conform discrete space $\mathcal{P}_{h}^{k} \subset \mathcal{P}$ as
\begin{equation}
\label{eq:qh1} \mathcal{P}_{h}^{k}=\mathcal{S}_{h}^{k}  \cap \mathcal{P},
\end{equation}
we present the Galerkin finite element approximation  to the {
Problem P}: {\it
  Find $p_h \in \mathcal{P}_{h}^{k}$ such that
  \begin{equation}
    \label{eq:pvh}
    a(p_h,q_h)= f(q_h)
    \quad \forall \; q_h \in \mathcal{P}_{h}^{k}.
  \end{equation}
}

 As $ \mathcal{P}_{h}^{k} \subset \mathcal{P}$, this problem has existence and uniqueness
of solution guaranteed from Lax-Milgram Lemma and the following
estimate holds \citep[see, for example, in][]{CIARLET78}
\begin{equation}
  \label{eq:errop}
  \| p - p_h \| + h \| \nabla p - \nabla p_h \| \leq C h^{k+1}
  \left| p \right|_{k+1} \, .
\end{equation}

After solving (\ref{eq:pvh}), we can evaluate the velocity field
directly through Darcy's law
\begin{equation}
  \label{eq:vg}
  \bm{u}_G = -\Kappa \nabla p_h.
\end{equation}
From (\ref{eq:errop}) the following error bound applies to this
approximation
\begin{equation}
  \label{eq:errovg}
  \| \bm{u} - \bm{u}_G \| \leq  Ch^k | p|_{k+1}.
\end{equation}
This direct approximation for the velocity field is, in principle,
completely discontinuous on the interface of the elements and
presents very poor mass conservation.

%
%
\subsection{Stabilized Mixed Formulations} \label{sec:stable}
%
%

A usual way to compute more accurate velocity fields is to use mixed
formulations, where the velocity and the potential fields are
approximated simultaneously. The lack of stability of the classical
dual mixed formulation \citep{RAVIART77} in  Lagrangian finite
element spaces \citep{FORTIN91} is one of the main reasons for the
development of Stabilized Mixed Formulations. Here we review some of
these stabilizations: the Galerkin Least-Squares (GLS) formulation
equivalent to a minimization introduced by \citet{LOULA88}
(Minimization/Galerkin Least-Squares Method - MGLS), the adjoint
stabilization proposed by \citet{MASUD2002} (Hughes Variational
Method for Darcy Flow - HVM) and the unconditionally stable method
CGLS presented in \cite{CORREA2008A}. All these formulations  are
free of mesh dependent parameters of stabilization and allow the use
of continuous Lagrangian elements for potential and velocity,
including same order interpolation for both variables. For a
comparative study of such methods, see \cite{CORREA2008A}.
As in this kind of mixed formulations the velocity is the primary
variable, we consider here only Neumann boundary conditions, that
are essential boundary conditions for the associated variational
principles.
%
%
%

Using $\mathcal{P}_{h}^{k}$ for the potential, as defined in (\ref{eq:qh1})  with
${\mathcal{P}}=H^1(\Omega)\backslash \mathbb{R}$ and defining the
$C^0(\Omega)$ Lagrangian finite element space for the velocity
\begin{equation}
\mathcal{U}_{h}^{l} =[\mathcal{S}_{h}^{l}]^2\cap \mathcal{U},
\end{equation}
where $\mathcal{U}={H_0(\div)}$ for MGLS and HVM, and $\mathcal{U}={H_0(\div)}\cap H(\curl_\Lambda,\Omega)$ for CGLS,
with
\[
H(\curl_\Lambda,\Omega) =  \left\{\bm{u} \in \left[ L^2(\Omega) \right]^2; \; {\bm{\curl}}
(\Lambda \bm{u}) \in L^2(\Omega)\right\},
\]
 we can present all these methods through the following
general form:

\noindent {\it Find $\{\bm{u}_h,p_h\} \in \mathcal{U}_{h}^{l}\times\mathcal{P}_{h}^{k}$ such that
\begin{equation}
B_{\alpha}\left( \{\bm{u}_h,p_h\}  ;  \{\bm{v}_h,q_h\} \right)= F_{\alpha}\left(\{\bm{v}_h,q_h\} \right)
\quad \forall \;\{\bm{v}_h,q_h\} \in \mathcal{U}_{h}^{l}\times\mathcal{P}_{h}^{k}
\end{equation}
where
\begin{eqnarray}
\nonumber B_{\alpha}\left( \{\bm{u}_h,p_h\}  ;  \{\bm{v}_h,q_h\} \right)&=&
(\Lambda\bm{u}_h,\bm{v}_h)- (\Div \bm{v}_h, p_h) - \delta_0(\Div \bm{u}_h, q_h) \\
\nonumber && + \delta_1\left( K\left( \Lambda \bm{u}_h + \nabla
p_h\right), \delta_0\Lambda \bm{v}_h + \nabla q_h  \right) \\
\nonumber
&& + \delta_2\left(\lambda\Div \bm{u}_h,\Div \bm{v}_h \right) \\
\label{eq:g-bilinear} && +\delta_3\left(K\curl(\Lambda \bm{u}_h),\curl
(\Lambda \bm{v}_h) \right)
\end{eqnarray}
and
\begin{equation}
F_{\alpha}\left(\{\bm{v}_h,q_h\} \right)= -\delta_0(f,q_h)+\delta_2(\lambda f,
\Div \bm{v}_h),
\end{equation}
}
with $\Lambda = \Kappa^{-1}$, $K = |\Kappa|_{\infty}$ and $\lambda =
K^{-1}$. In this abstract problem, the index $\alpha$ represents the
methods. In Table \ref{tab:stab} we present the values of the
parameters $\delta_i$ that define each method and their respective
{\em a priori} orders of convergence for the errors of the potential,
velocity and divergence of velocity, measured in the $L^2(\Omega)$
norm, when same order interpolations are used ($l=k$), for
sufficiently regular solutions, i.e. for both $\bm{u}$ and $p$ $\in
H^{k+1}(\Omega)$.

\vspace{.2cm}
\begin{center}
\begin{table}[htb]
\caption{Stabilized Mixed Methods} 
%
%
{\setstretch{1.35}
\begin{tabular}{  c | c|  c | c | c | c | c | c  }
\hline \hline
 $\alpha$ & $\delta_0$ & $\delta_1$ & $\delta_2$ & $\delta_3$ & $\|p-p_h\|$ & $\|\bm{u}-\bm{u}_h\|$ & $\|\Div \bm{u}-\Div \bm{u}_h\|$\\
\hline
MGLS  & $+1 $ & $\delta_1>0$ & $\delta_2>0$ &   0   & $O\left(h^{k+1}\right)$& $O\left(h^{k}\right)  $ & $O\left(h^{k}\right)  $ \\
HVM   & $-1$ & $+1/2$ & $ 0 $ &   0   & $O\left(h^{k+1}\right)$& $O\left(h^{k}\right)  $ & $O\left(h^{k-1}\right)$ \\
CGLS  & $+1$ & $-1/2$ & $1/2$ & $1/2$ & $O\left(h^{k+1}\right)$& $O\left(h^{k+1}\right)$ & $O\left(h^{k}\right)  $ \\
\hline\hline
\end{tabular}
}
\label{tab:stab}
%
%
\end{table}
\end{center}

%
%
\subsection{A Convergence Test} \label{sec:numerical1}
%
%

The orders of convergence presented in Table \ref{tab:stab} are
numerically confirmed in \cite{CORREA2008A} for both homogeneous and
nonhomogeneous media with regular solutions. We now test the
convergence of these methods when applied to a heterogeneous
and anisotropic media in which the conductivity tensor is
anisotropic and presents a smooth interface of discontinuity. In
this case the solution is piecewise regular, with a discontinuity in
the velocity field, which is the primary variable. In this study, we
adopted the parameters $\delta_1=\delta_2=1/2$ for the MGLS method.

This test problem from \cite{CRUMPTON95} is defined on the square
$[-1,1]\times[-1,1]$, with Dirichlet boundary conditions. The
conductivity is given by
  \[
 \Kappa= \left(
      \begin{array}{ccc}
    1 & \ \  0  \\
    0 & \ \ 1
      \end{array}
      \right)   \quad  x<0,  \qquad
   \Kappa=\gamma\left(
      \begin{array}{ccc}
    2 & \ \ 1  \\
    1 & \ \ 2
      \end{array}
      \right)   \quad x>0,
  \]
where the parameter $\gamma$ is used to vary the strength of the
discontinuity at $x=0$. In our experiments we used $\gamma=1.0$. The
exact potential field is given by
  \[
  p=
  \left\{
  \begin{array}{lcr}
  \gamma [2\sin(y)+\cos(y)] x + \sin(y), & & x<0 \\ \\
   \exp(x)\sin(y),                       & & x>0
  \end{array}
  \right.
  \]
and the velocity field can be directly calculated by Darcy's law
(\ref{eq:darcy1}) and (\ref{eq:darcy2}).
In the numerical study we prescribed only the normal component of
the velocity on the boundary.
The finite element solutions were computed by adopting uniform
meshes of $8 \times 8$, $16 \times 16$, $32 \times 32$ and $64
\times 64$ bilinear quadrilaterals (Q1) and $4\times4$, $8\times8$,
$16\times16$ and $32\times32$ biquadratic quadrilaterals (Q2). To
have sufficient accuracy in the numerical computation of the errors,
the numerical integration was performed with $3\times 3$ Gauss
quadrature for Q1 and $4\times4$ for Q2 elements.

\begin{figure}[htb]
\centering
\subfloat[Exact.]{\includegraphics[angle=0,width=.49\textwidth,angle=0]{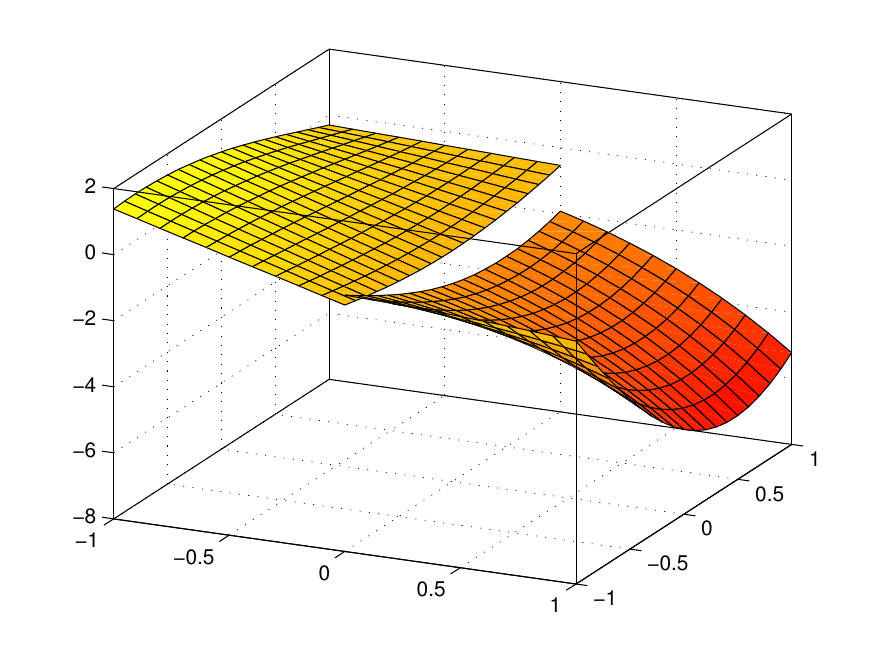}
\label{fig:vy-anis-exata}}
%
%
\subfloat[HVM.]{\includegraphics[angle=0,width=.49\textwidth,angle=0]{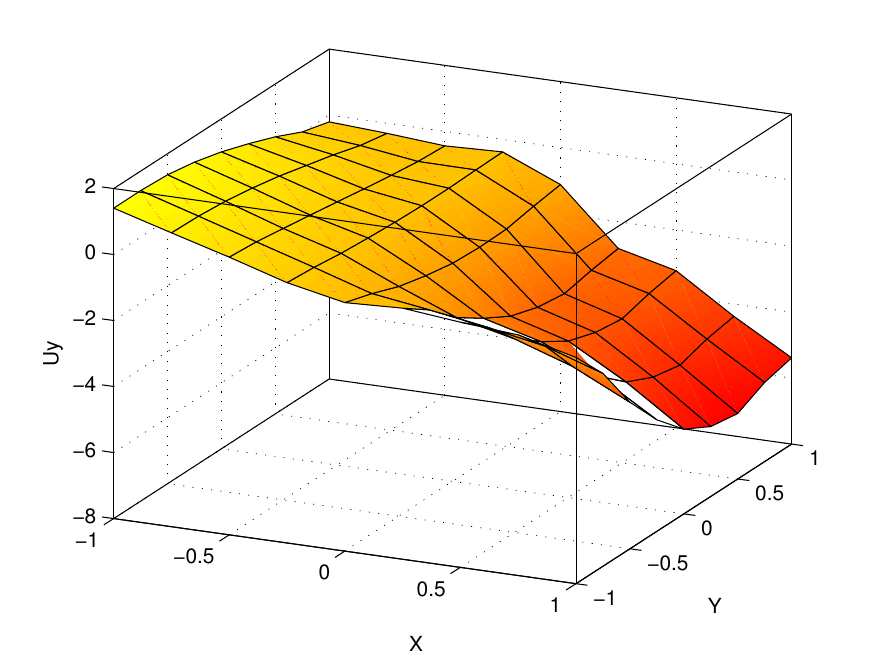}
    \label{fig:vy-anis-hvm}}
\caption{Component $u_y$;
    (a) Exact and (b) approximated by HVM with $8\times8$ bilinear elements.}
\end{figure}

The point in this example is that the discontinuity in the
conductivity leads to a velocity field in which the $y-$component is
discontinuous on the interface. Thus, the exact solution is not
sufficiently regular in the whole domain $\Omega$, but regular
enough in each subdomain $\Omega_1$ and $\Omega_2$.
As all stabilized mixed methods considered are based on
$C^0(\Omega)$ Lagrangian interpolation, we will investigate the
influence of imposing continuity on the approximate solution on the
interface $\Gamma$.
To illustrate this fact, the exact $y$ component of the velocity
field and its approximation by HVM with $8 \times 8$ bilinear
quadrilaterals are shown in Figures \ref{fig:vy-anis-exata} and
\ref{fig:vy-anis-hvm}, respectively.

\begin{figure}[htb]
\centering
\subfloat[]{\includegraphics[angle=0,width=.49\textwidth,angle=0]{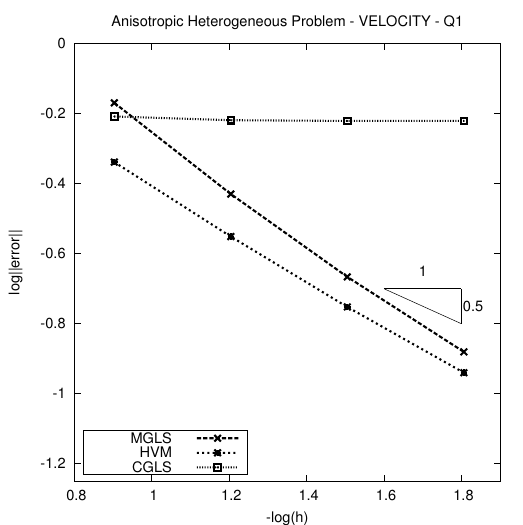}
\label{fig:vel-q1-continuo}}
\subfloat[]{\includegraphics[angle=0,width=.49\textwidth,angle=0]{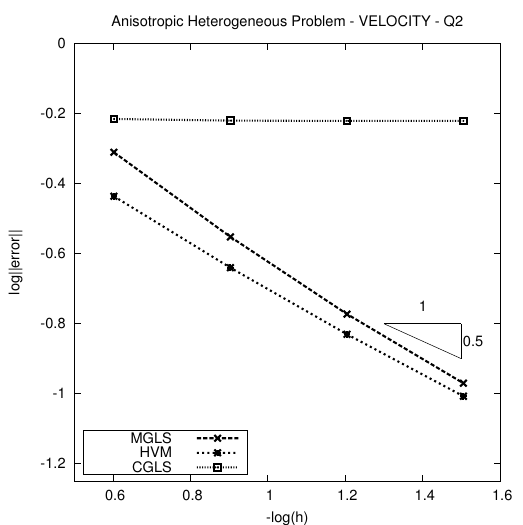}
\label{fig:vel-q2-continuo}}
\caption{Convergence study for an anisotropic heterogeneous problem
and continuous approximations;
Velocity with (a) bilinear elements and (b) biquadratic elements.}
\label{fig:vel-continuo}%
\end{figure}

In the convergence graphics, we present plots of the $\log$ of
$L^2(\Omega)$ norm of the errors {\em versus} $-\log(h)$.
Convergence results of $C^0(\Omega)$ velocity approximations for Q1
and Q2 elements are presented in Figures \ref{fig:vel-q1-continuo}
and \ref{fig:vel-q2-continuo}, respectively. Order of convergence
close to $O(h^{0.5})$ are obtained for HVM and MGLS methods for both
Q1 and Q2 elements, while no convergence is observed for the CGLS
method as a consequence of the lack of global regularity of the
exact solution combined with the use of continuous interpolations to
approximate a discontinuous field.

\begin{figure}[htb]
\centering
\subfloat[]{\includegraphics[angle=0,width=.49\textwidth,angle=0]{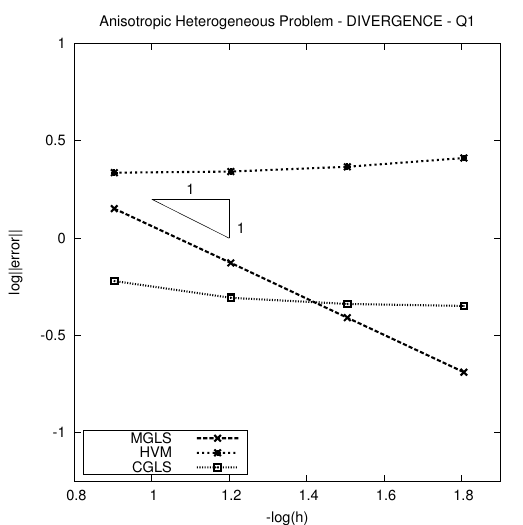}
\label{fig:div-q1-continuo}}
\subfloat[]{\includegraphics[angle=0,width=.49\textwidth,angle=0]{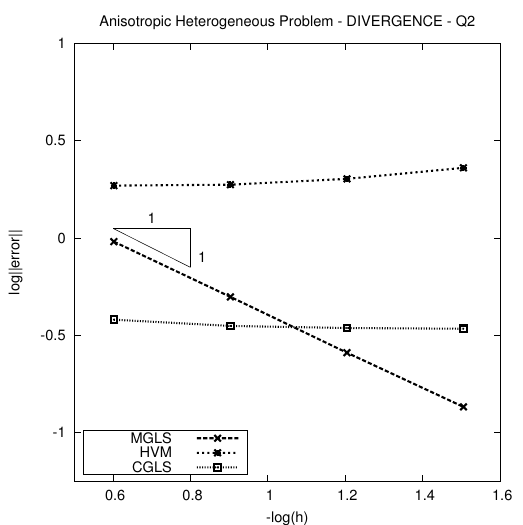}
\label{fig:div-q2-continuo}}
\caption{Convergence study for an anisotropic heterogeneous problem
and continuous approximations;
Divergence of velocity with (a) bilinear elements and (b) biquadratic elements.}
\label{fig:div-vel-continuo}%
\end{figure}

\begin{figure}[htb]
\centering
\subfloat[]{\includegraphics[angle=0,width=.49\textwidth,angle=0]{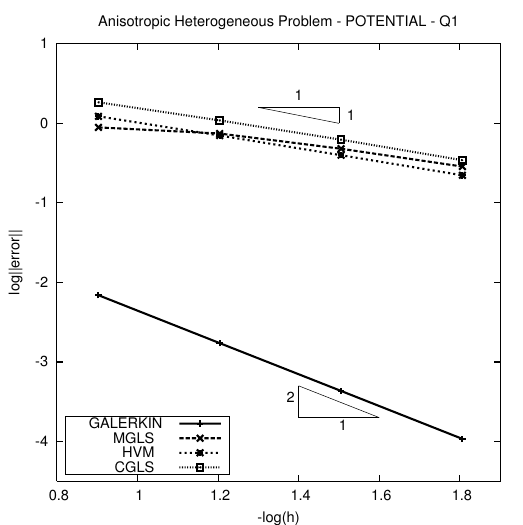}
\label{fig:pot-q1-continuo}}
\subfloat[]{\includegraphics[angle=0,width=.49\textwidth,angle=0]{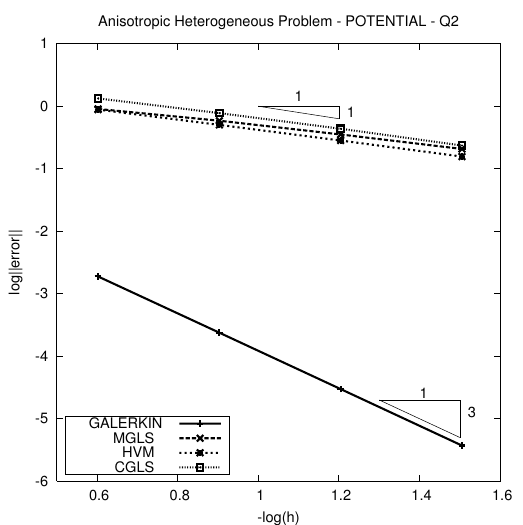}
\label{fig:pot-q2-continuo}}
\caption{Convergence study for an anisotropic heterogeneous problem
and continuous approximations;
Potential with (a) bilinear elements and (b) biquadratic elements.}
\label{fig:pot-continuo}%
\end{figure}

Figures \ref{fig:div-q1-continuo} and \ref{fig:div-q2-continuo}
present convergence results for the divergence of the velocity
approximations.
The poor accuracy of the approximations for the velocity on the
interface degrades the rates of convergence of the divergence for
all methods. Convergence rates close to $O(h^{1.0})$ are obtained
for MGLS with Q1 and Q2 elements.
Concerning HVM approximations, according to the predicted estimate
of convergence, we should not expect convergence for divergence
using Q1 elements even for regular exact solutions, except for the
``superconvergence'' usually observed with linear and bilinear $C^0$
velocity approximations.
But $O(h)$ convergence is expected with Q2 elements for regular
solutions. However, no convergence is observed. This can be
explicated by the error estimate for HVM (Table \ref{tab:stab})
combined with the lack of global regularity of the exact solution of
this problem and the use of $C^0(\Omega)$ Lagrangian interpolation
for the velocity field.

In Figures \ref{fig:pot-q1-continuo} and \ref{fig:pot-q2-continuo}
we plot the errors for potential using Q1 and Q2 elements,
respectively. The  methods show convergence rates of lower order
than those we would expected from Table \ref{tab:stab}, with errors
considerably large, when compared with the results of the single
field formulation (Galerkin).
This can also be explained by the coupled estimates for potential
and velocity approximations, characteristic of mixed methods,  and
the lack of global regularity of the exact solution.

To recover stability, convergence and accuracy of the stabilized
mixed methods presented in this section, when applied to our model
problem, we will introduce in the Section \ref{sec:discont} a
generalization of these formulations in which
continuity/discontinuty of the velocity and potential fields on the
interface $\Gamma$ are exactly imposed.

%
%
\section{Broken Space Formulations}
\label{sec:broken}
%
%

In this Section we comment on special single field formulations of
our model problem, based on a extended minimization problem and on a
discontinuous Galerkin formulation in which the continuity or
discontinuity of the velocity field are imposed in a weak sense.

\subsection{ Extended Minimization Problem}

We note that {Problem P} is equivalent to the following minimization
problem in $H^1_0(\Omega)$

\vspace{.25cm}
\paragraph{Problem J}{\it
  Find $p \in {\mathcal Q}=H^1_0(\Omega)$ such that

  \begin{equation}
    J(p) \le J(q)
    \quad \forall \quad  q \in {\mathcal{Q}}
  \end{equation}
  with
   \begin{equation}
 J(q) = \frac{1}{2} a(p,q) - f(q) .
   \end{equation}  %
}

 \cite{BEVI74} presented an equivalent
formulation for this minimization problem by constructing an
extended functional, defined on a broken function space, whose
minimum is also given by the solution of { Problem J}, with the
interface conditions naturally incorporated in the variational
formulation. The minimum of this extended functional can be
approximated by minimizing sequences of functions, discontinuous on
the interface $\Gamma$. To briefly present this formulation applied
to our model problem, we introduce some additional notations and
definitions.
 Let
\[
{\mathcal{V}}_i(\Omega_i) = \{ q\in H^2(\Omega_i) \, : \ q = 0 \
\hbox{on } \ \Gamma_i \, , \ i=1,2 \}
\]
and
\[
{\mathcal{V}}_e(\Omega) = \{q\in L^2(\Omega) ; q_i = q|_{\Omega_i} \in
{\mathcal{V}}_i(\Omega_i)\, , \ i=1,2 \}.
\]

Defining
\begin{multline}
J_e(q) = \sum_{i=1}^2\displaystyle \int_{\Omega_i} \left( \frac{1}{2}
\Kappa_i\nabla q_i \cdot \nabla q_i - f_iq_i \right) \mathrm{d}{\Omega} \\
- \frac{1}{2} \int_{\Gamma} (\Kappa_2 \nabla q_2 + \Kappa_1 \nabla
q_1)\cdot \bm{n} (q_2 - q_1) \mathrm{d}s
\end{multline}
the following minimization problem is set
\vspace{.25cm}
\paragraph{Problem Je} {\it
  Find $p \in {\mathcal{V}}_e(\Omega) $ such that

  \begin{equation}
    J_e(p) \le J_e(q)
    \quad \forall \quad  q \in {\mathcal{V}}_e(\Omega).
  \end{equation}
}
The minimum condition of {Problem Je} leads to
\begin{multline*}
\sum_{i=1}^2\int_{\Omega_i} [-\Div (\Kappa_i\nabla p_i) - f_i]q_i
\mathrm{d}{\Omega} 
+\frac{1}{2} \int_{\Gamma} (\Kappa_2 \nabla p_2 - \Kappa_1 \nabla
p_1)\cdot \bm{n} (q_2 + q_1) \mathrm{d}s
\\
-\frac{1}{2} \int_{\Gamma} (\Kappa_2 \nabla q_2 + \Kappa_1 \nabla
q_1)\cdot \bm{n} (p_2 - p_1) \mathrm{d}s =0
\end{multline*}
and, since the variations are arbitrary, we have:
\begin{equation}
-\Div (\Kappa_i\nabla p_i) = f_i \, , \quad i=1,2 \, , \quad \hbox{in}
\quad \Omega
\end{equation}
\begin{equation}
\label{eq:cond1} (\Kappa_2 \nabla p_2 - \Kappa_1 \nabla p_1)\cdot \bm{n} =0
\, , \quad \hbox{on} \quad \Gamma
\end{equation}
\begin{equation}
p_2 - p_1 =0 \, , \quad \hbox{on} \quad \Gamma.
\end{equation}
We observe that condition (\ref{eq:cond1}) is identical to condition
(\ref{eq:nflx}) posed in terms of the potential, by using Darcy's
law.

In   \cite{BEVI74}  this formulation is applied to the torsion
problem of a rod composed of two different materials using Ritz's
method with a minimizing sequence whose coordinate functions are
discontinuous on the interface of the materials. In their
conclusions the authors say:

{\it ``... the use of such functionals to derive the fundamental
equations of a finite element technique is promising, in those case
where it is convenient to relax the requirement of continuity of
some of the variables''}.

We may  naturally associate this conclusion with the basic idea of
Discontinuous Galerkin Finite Element Methods \citep[e.
g.,][]{DDUPONT76,WHEELER78,ARNOLD82,ODEN98,DUTRA02,ARNOLD02}.
Certainly, \cite{BEVI74} did not try to implement their promising
idea in the IBM 1130 with 32K of RAM memory they had access in
COPPE/UFRJ.

In the next section we present a weak formulation of our model
problem, in broken function spaces on a typical subdivision of the
domain $\Omega$, based on which the hp discontinuous Galerkin finite
element method is constructed.

\subsection {Weak Formulations in Broken Spaces}

Let $\mathcal{T}$ be a subdivision of $\Omega$ into $Ne$ elements $\kappa$
such that the interface $\Gamma$ is the union of elements edges $e$.
Denote by $\mathcal{E}$ the set of all edges associated with the subdivision
$\mathcal{T}$ and define $\mathcal{E}_{\partial} = \partial \Omega$ and $\mathcal{E}_{int}$
as the set of boundary and interior edges of the subdivision $\mathcal{T}$,
respectively. For a given edge $e \in \mathcal{E}_{int}$ adjacent to the
elements $\kappa=\kappa_i$ and $\kappa^{\prime}=\kappa_j$, $i>j$, we
define
$$
\jp{p}_e = p|_{\partial\kappa\cap e} - p|_{\partial\kappa^{\prime}\cap
e}
$$
$$
\av{p}_e = \frac{1}{2}\left (p|_{\partial\kappa\cap e}
+ p|_{\partial\kappa^{\prime}\cap e}\right )
$$
where $\partial \kappa$ denotes the union of all edges of the
element $\kappa$.

 The broken weak formulation of our model problem reads
\vspace{.2cm}
\paragraph{Problem B} {\it
  Find $p \in \mathcal{P}^*(\Omega, \mathcal{T}) $ such that

  \begin{equation}
    b(p,q)= l(q)
    \quad \forall \quad  q \in H^2(\Omega,\mathcal{T})
  \end{equation}
}
  with
  \begin{equation}
    H^2(\Omega,\mathcal{T})=\left\{ q\in L^2(\Omega) ;
    \left. q\right|_{\kappa}\in H^2(\kappa) \;\forall\; \kappa \in \mathcal{T}
    \right\} \, ,
  \end{equation}
\begin{equation}
\mathcal{P}^*(\Omega, \mathcal{T}) = \{ q\in H^2(\Omega,\mathcal{T})\, : q, \Kappa \nabla q\cdot \bm{n} \hbox{ are continuous on each edge} \ e \in \mathcal{E}_{int} \} \, ,
\nonumber
\end{equation}
\begin{equation}
b(p,q)= \sum_{\kappa\in \mathcal{T}} \int_{\kappa} \Kappa \nabla p \cdot \nabla
q \,
\mathrm{d}{\Omega} + b_{\Gamma} (p,q) + b_{\partial} (p,q) \, ,
\end{equation}
\begin{equation}
l(q)= \sum_{\kappa\in \mathcal{T}} \int_{\kappa} fq \mathrm{d}\Omega
   \end{equation}
   where
 \begin{equation}
  b_{\Gamma}(p,q)=\int_{\Gamma_{int}} \left(\alpha \jp{p} \av{\Kappa\nabla  q \cdot \bm{n}} - \av{\Kappa\nabla  p \cdot \bm{n}}\jp{q} + \beta \jp{p}\jp{q}\right) \mathrm{d}s
  \end{equation}
 \begin{equation}
  b_{\partial}(p,q)= \int_{\partial\Omega} \left(\alpha p \av{\Kappa\nabla  q \cdot \bm{n}} - \av{\Kappa\nabla  p \cdot \bm{n}} q + \beta pq\right)
  \mathrm{d}s \, .
 \end{equation}
 The bilinear form  $b_{\Gamma} (p,q)$
 represents the jump terms on the interior edges of the elements
 $\Gamma_{int}$ responsible for weak continuity across the element
 edges, while  $b_{\partial}(p,q)$, defined on the boundary
 $\partial\Omega$, enforces Dirichlet boundary condition in a weak
 sense.

 The free parameters $\alpha$ and $\beta$ play an important role in the
 stability of discontinuous formulations; see \cite{ODEN98,DUTRA02,ALVAREZ2006}.
 We should note that $\beta$ is a penalty parameter which is usually of
 the order of $h^{-1}$ while $\alpha$ is normally considered in the
 interval $ [ -1, 1]$ \citep{ODEN98,DUTRA02,ARNOLD02}.

The choice $\alpha =1$, though leading to a nonsymmetric bilinear
form, is usually adopted for originating a more stable method than
the symmetric one. For this choice, we have
   \begin{equation}
  b(q,q)= \sum_{\kappa\in \mathcal{T}} \int_{\kappa} \Kappa \nabla q \cdot \nabla
  q \,
  \mathrm{d}{\Omega}+ \int_{\Gamma_{int}} \beta \jp{q}^2 \mathrm{d}s
  + \int_{\partial\Omega}  \beta q^2 \mathrm{d}s
  \end{equation}
which proves the stability of this formulation in the broken space.

 For $\alpha=-1$ the bilinear form $ b(\cdot, \cdot)$ preserves
 symmetry and {Problem B} becomes equivalent to the minimization
 problem associated with the functional
  \begin{equation}
   J_b(q) = \frac {1}{2} b(q,q) - l(q) \, ,
   \end{equation}
in the broken space $H^2(\Omega,\mathcal{T})$. Stability of this symmetric
formulation is also proved for appropriate choice of the penalty
parameter $\beta$ \citep{ARNOLD02}.

Defining $H^2_0(\Omega,\mathcal{T})=\{q\in H^2(\Omega,\mathcal{T}), q=0 \mbox{ on }
\partial\Omega \} $ and taking $\beta= 0$, we present a natural
generalization of the extended minimization problem of \cite{BEVI74}
as follows:
\vspace{.2cm}
\paragraph{Problem Jb} {\it
  Find $p \in H^2_0(\Omega,\mathcal{T}) $ such that

  \begin{equation}
    J_b(p) \le J_b(q)
    \quad \forall \quad  q \in H^2_0(\Omega,\mathcal{T})
  \end{equation}
}
with
\begin{equation}
  J_b(q) = \frac{1}{2} \sum_{\kappa\in \mathcal{T}} \int_{\kappa} \Kappa \nabla q \cdot \nabla
  q \,
  \mathrm{d}{\Omega}- \int_{\Gamma_{int}} \av{\Kappa\nabla  q \cdot \bm{n}}\jp{q} \mathrm{d}s
  - \sum_{\kappa\in \mathcal{T}} \int_{\kappa} fq \mathrm{d}{\Omega}\, .
\end{equation}
This minimization problem is equivalent to Problem B with
$\alpha=-1$ and $\beta=0$, restricted to $H^2_0(\Omega,\mathcal{T})$, due to
the fact that in \cite{BEVI74} Dirichlet boundary condition is
imposed in the strong sense.

 Recently, a renewed interest in discontinuous Galerkin methods has surged for their
flexibility in hp adaptivity strategies, naturally allowing
variation of the polynomial degrees of the interpolation functions
and irregular meshes since no pointwise continuity is imposed at the
element interfaces. Additionally, discontinuous Galerkin methods
have very interesting local properties, such as conservation, and
additional stability compared to classical continuous Galerkin
methods \citep{DUTRA02,ALVAREZ2006,BREZZI2006}. However,
discontinuous Galerkin methods present a serious drawback compared
to their continuous counterparts. The explosion of degrees of
freedom. This has motivated the selective use of the discontinuous
Galerkin formulation given rise to continuous/discontinuos
formulations, as in \cite{DAWSON2005,DUARTE2000} and
\cite{CORREA2007}.

%
%
\section{Continuous/Discontinuous Formulations}
\label{sec:discont}
%
%

As we demonstrated in Section \ref{sec:stable}, the use of
$C^0(\Omega)$ Lagrangian interpolations is not adequate for problems
in which the primal variable is discontinuous. In Section
\ref{sec:broken} we discussed some alternative formulations to deal
with this kind of problem, using discontinuous basis functions
associated with extended variational principles.
In this Section we propose a generalization of the $C^0(\Omega)$
Lagrangian based stabilized mixed formulations in which continuity
or discontinuity of the velocity and potential fields on the
interface $\Gamma$ are exactly imposed, recovering stability,
accuracy and the optimal rates of convergence of classical
Lagrangian interpolations.

%
%
\subsection{Constrained Variational Formulation}
\label{sec:const}
%
%

%
For simplicity we limit this presentation to the CGLS stabilization
with velocity and potential regular enough in each subdomain
$\Omega_i$, $i=1,2$. Numerical results confirm that, with this new
formulation, we obtain rates of convergence equivalent to the
corresponding stabilizations of the homogeneous problem with regular
solution.
Let
\[
 \mathcal{U}_j =H_{0,\Gamma_j}(\Div,\Omega_j)
 \cap H( \curl\Lambda_j,\Omega_j)
 \]
 with
\[
H_{0,\Gamma_j}(\Div,\Omega_j)= \left\{ \bm{u} \in H(\Div,\Omega_j);
\bm{u}\cdot\bm{n}=0 \mbox{ on } \Gamma_j \right\},
\]
\[
H(\curl \Lambda_j ,\Omega_j)= \left\{\bm{u} \in \left[ L^2(\Omega_j)
\right]^2; \; \curl (\Lambda_j \bm{u}) \in L^2(\Omega_j)\right\},
\]
and
\[
\mathcal{P}_j= H^1(\Omega_j)
\]
be the velocity and potential spaces in the subdomain $\Omega_j$,
$j=1,2$, respectively. Then, defining
\[
{\mathcal S}_{hj}^l=\left\{ \varphi_h \in C^0(\Omega_j)\;;\;
\left.\varphi_h \right|_{\Omega^e}\in \mathbb{P}_{l,l}(\Omega^e) \;
\forall\; \Omega^e \in {\mathcal T}_h \cap \Omega_j \right\}
\]
\[
\mathcal{U}^l_{hj}= \left[{\mathcal S}_{hj}^l\right]^2 \cap \mathcal{U}_j
\]
\[
\mathcal{P}^k_{hj}= {\mathcal S}_{hj}^k \cap \mathcal{P}_j
\]
we have the following finite element spaces for the admissible
variables
\[
\mathcal{U}_{h}^{l}=\left\{\bm{v}_h\; ; \; \bm{v}_h^1=\left.\bm{v}_h\right|_{\Omega_1} \in
\mathcal{U}_{h1} \;,\; \bm{v}_h^2=\left.\bm{v}_h\right|_{\Omega_2} \in \mathcal{U}_{h2} \right\}
\]
\[
\mathcal{P}_{h}^{k}=\left\{ q_h\; ; \; q_h^1=\left.q_h\right|_{\Omega_1} \in
\mathcal{P}^k_{h1} \;,\; q_h^2=\left.q_h\right|_{\Omega_2} \in \mathcal{P}^k_{h2}
\right \} .
\]
We will also use the spaces of continuous functions
\[
\mathcal{U}^l_{h0}=\mathcal{U}_{h}^{l} \cap \left[ C^0(\Omega)\right]^2\; ,\qquad
\mathcal{P}^k_{h0}=\mathcal{P}_{h}^{k}\cap C^0(\Omega).
\]

In these spaces, a natural formulation of the CGLS stabilization
reads as
\vspace{.25cm}
\paragraph{Problem CGLS} {\it Find $\{\bm{u}_h,p_h\}
\in \mathcal{U}_{h}^{l}\times \mathcal{P}_{h}^{k}$  such that $\forall \; \{\bm{v}_h,q_h\} \in \mathcal{U}^l_{h0} \times
\mathcal{P}^k_{h0}$
\begin{multline*}
B_\mathrm{CGLS}^1(\{\bm{u}^1_h,p^1_h\},\{\bm{v}_h,q_h\}) +
B_\mathrm{CGLS}^2(\{\bm{u}^2_h,p^2_h\},\{\bm{v}_h,q_h\}) = \\
F^1_\mathrm{CGLS}(\{\bm{v}_h,q_h\})  +F^2_\mathrm{CGLS}(\{\bm{v}_h,q_h\}) 
\end{multline*}
with
\begin{eqnarray*}
B^i_\mathrm{CGLS}(\{\bm{u}_h,p_h\},\{\bm{v}_h,q_h\}) =   (\Lambda \bm{u}_h,\bm{v}_h)_{\Omega_i}
- (\Div\bm{v}_h,p_h)_{\Omega_i} - (\Div\bm{u}_h,q_h)_{\Omega_i}  \\
-\frac{1}{2} \left( K\left( \Lambda \bm{u}_h + \nabla p_h\right),
\Lambda \bm{v}_h + \nabla q_h\right)_{\Omega_i}\\
+\frac{1}{2} (\lambda\Div \bm{u}_h ,\Div \bm{v}_h)_{\Omega_i} +
\frac{1}{2}(K\curl (\Lambda \bm{u}_h), \curl (\Lambda
\bm{v}_h))_{\Omega_i} \,,
\end{eqnarray*}
\[
F^i_\mathrm{CGLS}(\{\bm{v}_h,q_h\})= -(f,q_h)_{\Omega_i}+\frac{1}{2}(\lambda f,
\Div \bm{v}_h)_{\Omega_i}\;,
\]
plus the discrete interface conditions
\begin{equation}
  \label{eq:restr-p}
  p^1_h = p^2_h
  \quad \mbox{ on } \quad \Gamma
\end{equation}
\begin{equation}
  \bm{u}^1_h\cdot \bm{n} = \bm{u}^2_h\cdot \bm{n} \quad \mbox{ on } \quad
  \Gamma
  \label{eq:restr-un}
\end{equation}
\begin{equation}
  \Lambda_1\bm{u}^1_h\cdot \tau= \Lambda_2  \bm{u}^2_h\cdot \tau
  \quad \mbox{ on } \quad \Gamma.
  \label{eq:restr-ut}
\end{equation}
}

One typical way to impose the above conditions on the interface
$\Gamma$ is using Lagrange multiplier. We propose a simpler method
consisting in imposing the interface conditions at element level
leading to a global problem similar to the corresponding
$C^0(\Omega)$ finite element formulations presented in Section
\ref{sec:stable}.

\subsection{ An Auxiliary Unconstrained Formulation}
\label{sec:unconstrained}

Let $\bm{n}_1=\bm{n}$ be the outward normal to the boundary of $\Omega_1$
at a generic point  $P(\bm{x})$ on the interface $\Gamma$. Thus, for
$\epsilon>0$
\[
\bm{u}^1 = \lim_{\epsilon \rightarrow 0}\bm{u}(\bm{x}-\epsilon \bm{n})
\]
\[
\bm{u}^2 = \lim_{\epsilon \rightarrow 0}\bm{u}(\bm{x}+\epsilon \bm{n}) \, .
\]
We should observe that
\[
 \lim_{\epsilon \rightarrow 0} \Kappa^{-1}(\bm{x}-\epsilon \bm{n})\bm{u}(\bm{x}-\epsilon \bm{n})\cdot\tau=
                               \Kappa_1^{-1}\bm{u}^1 \cdot \tau = \Kappa_2^{-1}\bm{u}^2 \cdot \tau =
 \lim_{\epsilon \rightarrow 0} \Kappa^{-1}(\bm{x}+\epsilon \bm{n})\bm{u}(\bm{x}+\epsilon \bm{n})\cdot \tau
\]
which implies \( \displaystyle \bm{u}^1 \cdot \tau \ne \bm{u}^2 \cdot \tau \) for
$\Kappa_1 \ne \Kappa_2 $.  Thus, the $C^0(\Omega)$ Lagrangian
interpolation for the velocity field is obviously incompatible with
this last condition since it imposes exactly
\begin{equation}
  \label{eq:resc}
  \bm{u}^1_h\cdot \bm{n} = \bm{u}^2_h\cdot \bm{n} \quad \mbox{ and } \quad
  \bm{u}^1_h\cdot \tau= \bm{u}^2_h\cdot \tau
\end{equation}
on the interfaces of the elements.
In this case, the approximate solution ${\bm{u}}_h$ at a global node on
$\Gamma$ will represent an intermediate value of the discontinuous
solution. Our idea is to impose exactly the constraints
(\ref{eq:restr-p}), (\ref{eq:restr-un}) and (\ref{eq:restr-ut}) on
the interface of discontinuity and assemble the global finite
element approximation in terms of a continuous reference
``solution'' $ \{ \bar{{\bf u}}_h , \bar {p}_h \}$ related with the
true solution $ \{ {{\bf u}}_h , {p}_h \}$  through a linear
transformation. This approach leads to a finite element formulation
with the same connectivity of the stabilized mixed $C^0(\Omega)$
formulation in which the reference pair $ \{ \bar{{\bf u}}_h , \bar
{p}_h \}$ is obtained solving the following problem:

\vspace{.2cm}
\paragraph{Problem $\mathrm{CGLS}_{\pi}$} {\it Find $\{ \bar{\bm{u}}_h,\bar{p}_h\}
\in \mathcal{U}^l_{h0} \times \mathcal{P}^k_{h0}$  such that
\begin{eqnarray*}
B^{\pi}_\mathrm{CGLS}(\{\bar{\bm{u}}_h,\bar{p}_h\},\{\bar{\bm{v}}_h,\bar{q}_h\})=
F^{\pi}_\mathrm{CGLS}(\{\bar{\bm{v}}_h,\bar{q}_h\}) \quad \forall \;
\{\bar{\bm{v}}_h, \bar{q}_h \} \in \mathcal{U}^l_{h0} \times \mathcal{P}^k_{h0}
\end{eqnarray*}
with
\begin{eqnarray*}
B^{\pi}_\mathrm{CGLS}(\{\bar{\bm{u}}_h,\bar{p}_h\},\{\bar{\bm{v}}_h,\bar{q}_h\})=
B_\mathrm{CGLS}^1(\{\pi\bm{u}^1_h,p_h\},\{\pi\bm{v}_h^1,q_h\})&+& \\
B_\mathrm{CGLS}^2(\{\pi\bm{u}^2_h,p_h\},\{\pi\bm{v}_h^2,q_h\}) &&
\end{eqnarray*}
\[
F^{\pi}_\mathrm{CGLS}(\{\bar{\bm{v}}_h,\bar{q}_h\})=
F^1_\mathrm{CGLS}(\{\pi\bm{v}_h^1,q_h\})+F^2_\mathrm{CGLS}(\{\pi\bm{v}_h^2,q_h\})
\]
and $\pi$ denoting a linear transformation from $\mathcal{U}^l_{h}$ to
$\mathcal{U}^l_{h0}$ defined by strongly imposing the interface conditions.
 }

In the next section, taking $\Omega_2$ as a reference subdomain,  we
define the linear transformation $\pi: \mathcal{U}^l_{h} \to \mathcal{U}^l_{h0}$ and
construct the element matrices and load vectors of the elements
belonging to $\Omega_1$ adjacent to the interface $\Gamma$, which
are the only elements whose matrices and vectors will be  affected
by this transformation.

\subsection{ Imposing the Interface Conditions}
\label{sec:idi}

We start by choosing the reference solution $\bar{\bm{u}}_h|_{\Omega_2}
= \bm{u}_h^2$ and $\bar{p}_h|_{\Omega_2} = p_h^2 $. Thus, in matrix
form we write
\begin{equation}
\label{eq:uref1}
 U^2_h=\left\{
\begin{array}{c}
\bm{u}^2_h \\
p^2_h
\end{array}
\right\} = \bar{U}_h= [I]\left\{
\begin{array}{c}
\bar{\bm{u}}_h \\
\bar{p}_h
\end{array}
\right\},
\end{equation}
where $[\mathrm{I}]$ is the identity. This equation relates $U_h^1$
to the reference solution. From equations (\ref{eq:restr-p}),
(\ref{eq:restr-un}) and (\ref{eq:restr-ut}) we can write
\[
\mathrm{T}_1 U_h^1= \mathrm{T}_2  \,U_h^2 =  \mathrm{T}_2
\,\bar{U}_h \quad \mbox{ on } \quad \Gamma,
\]
with the Cartesian matrix of $ \mathrm{T}_i $ given by
\begin{equation}
\label{eq:tmatrix}
 \left[ \mathrm{T}_i \right] = \left[
\begin{array}{ccccc}
  \lambda_{11}^i \tau_1 + \lambda_{21}^i\tau_2  & &
  \lambda_{12}^i \tau_1 + \lambda_{22}^i\tau_2  & & 0 \\
   \mathrm{n}_1 & & \mathrm{n}_2 & & 0 \\
   0 & & 0 & & 1
\end{array}
\right]
\end{equation}
where $(\tau_1,\tau_2)$ and $(\mathrm{n}_1,\mathrm{n}_2)$ are the
Cartesian components of the unit vectors $\tau$ and $\bm{n}$,
respectively, and $\lambda_{lm}^i$ are the components of the
hydraulic resistivity tensor $\Lambda_i=\Kappa_i^{-1}$. Using the fact
that  $\mathrm{T}_1$ (and also $\mathrm{T}_2$) is invertible, as a
consequence of
\[
\mathrm{det}\left[\mathrm{T}_i\right]=\Lambda_i \tau\cdot\tau>0,
\]
we get
\begin{equation}
\label{eq:uref2} U_h^1 = \mathrm{T}_1^{-1} \mathrm{T}_2\; \bar{U}_h
\quad \mbox{ on } \quad \Gamma.
\end{equation}

Introducing the  matrices
\begin{equation}
\label{eq:qmatrix}
 \left[ \mathrm{Q}_i \right] = \left[
\begin{array}{ccc}
  \lambda_{11}^i \tau_1 + \lambda_{21}^i\tau_2  & &
  \lambda_{12}^i \tau_1 + \lambda_{22}^i\tau_2   \\
   \mathrm{n}_1 & & \mathrm{n}_2
\end{array}
\right]
\end{equation}
for $i=1,2$, we can define the linear transformation $\pi$ at each
node of the mesh, and consequently at any point of the subdomains
$\Omega_1$ and $\Omega_2$ using the finite element interpolation
functions, as follows:
\begin{enumerate}
{\it
\item[1.] $[\pi] = [I]$ at any node of the mesh belonging to
$\Omega_2$ or $\Omega_1/\Gamma$,

\item[2.] $[\pi] = [Q] = [Q_1]^{-1} [Q_2] $ at any node of the mesh belonging to $\Omega_1\cap
\Gamma$.
}
\end{enumerate}

Using equations  (\ref{eq:uref1}) and  (\ref{eq:uref2}) we are able
to pose the interface problem in terms of a continuous global vector
of unknowns $\{\bar{U}\}$ only.  This can be done at element level
as follows. Let $\{U^e\}$ be the array of unknowns of a generic
element $\Omega^e$ with just one edge belonging to the interface of
discontinuity, and $ \left\{ U^e_{\Gamma} \right\}$ the unknowns
associated with the nodes on that edge belonging to the interface.
By using (\ref{eq:uref2}) to each node on $\Gamma$, we can construct
a matrix $\left[ \displaystyle T_\Gamma \right]$ and write
\begin{equation}
  \label{eq:rel1}
  \displaystyle \left\{ \displaystyle U^{e^{\,}}_{\Gamma} \right\} = \left[ \displaystyle T_\Gamma \right]
  \{ \displaystyle U^{e'}_{\Gamma} \}
\end{equation}
where $\Omega^e$ and $\Omega^{e'}$ are elements that share the
interface of material discontinuity. Choosing $\{U^{e'}_{\Gamma}\}$
as reference we have
\[
  \left\{ U^e \right\} = \left\{
  \begin{array}{c}
    \left\{ U^e_{\bar{\Omega} \backslash \Gamma} \right\} \\
    --- \\
    \left[ \displaystyle T_\Gamma \right] \{ \displaystyle U^{e'}_{\Gamma} \}
  \end{array}
  \right\}
     = \left[ T \right]
  \left\{
  \begin{array}{c}
    \left\{ U^e_{\bar{\Omega} \backslash \Gamma} \right\} \\
    --- \\
    \{ \displaystyle U^{e'}_{\Gamma} \}
  \end{array}
  \right\}
\]
or
\begin{equation}
  \label{eq:trans1}
  \left\{ U^e \right\} = \left[ T \right] \{ \bar{U}^{e} \}
\end{equation}
where  $[T]$ is the matrix of the linear transformation mapping  $\{
\bar{U}^{e} \}$ in $ \{ U^e \}$. The contribution of the element
$\Omega^e$ to the global stiffness $[K]$ is given by:
\[
[\tilde{K}^e]\{U^e\}=[ T]^{T} [ K^e ] [ T ] \{\bar{U}^e\}
            =[ \bar{K}^e ] \{\bar{U}^e\},  
\]
with the new matrix  $[\bar{K}^e] = [ T]^{T} [ K^e ] [ T ]$
preserving the symmetry. The corresponding load vector is $
\{\bar{F}^e \} = [ T]^{T} \{F^e\}$.

The continuous/discontinuous approximation $\{ U^e \}$ on the
interface is obtained at element-level  using (\ref{eq:trans1})
after solving the global system in  $\{\bar{U}\}$.

%
%
\subsection{Stability }
%
%

We observe that { Problem $\mathrm{CGLS}_{\pi}$} is a compatible
$C^0$ formulation that fits in the class of stabilized mixed methods
described in Section 3.2.  For the CGLS method, the following result
on stability in the sense of Babu\v ska is proved in
\cite{CORREA2008A}:
{\it There exists a constant $\alpha > 0 $ such that
$\forall \; \{{\bm{v}}_h, {q}_h \} \in
\mathcal{U}^l_{h0} \times \mathcal{P}^k_{h0}$
\begin{equation}
\label{eq:cp-coerc}
 \sup_{\{\bar{\bm{v}}_h, \bar{q}_h \}  \in \mathcal{U}^l_{h0} \times \mathcal{P}^k_{h0}}
 \frac{B^{\pi}_{\mathrm{CGLS}}(\{\bar{\bm{u}}_h,\bar{p}_h\},\{\bar{\bm{v}}_h,\bar{q}_h\})}{\| \{\bar{\bm{v}}_h, \bar{q}_h
 \}\|_{\mathcal{U}\times\mathcal{P}}}\geq  \alpha \| \{ \pi\bar{\bm{u}}_h,\bar{p}_h\} \|_{\mathcal{U}\times\mathcal{P}}
\end{equation}
where
\begin{equation}
\label{eq:cpf} \| \{\bm{u},p\} \|_{\mathcal{U}\times\mathcal{P}}s = \| \bm{u} \|^2 + \| \Div \bm{u} \|^2 + \|
\curl \Lambda \bm{u} \|^2 + \|\nabla p \|^2 \,.
\end{equation}
}
Considering that $\pi$ is continuous and bounded bellow, we
have
\begin{equation}
\label{eq:cpi}
 \| \{ \pi\bar{\bm{u}}_h,\bar{p}_h\} \|_{\mathcal{U}\times\mathcal{P}} \ge \bar{\beta} \|
\{ \bar{\bm{u}}_h,\bar{p}_h\} \|_{\mathcal{U}\times\mathcal{P}}
\end{equation}
with $\bar{\beta} > 0$, and the stability of { Problem
$\mathrm{CGLS}_{\pi}$} is proved. As the bilinear form
$B_\mathrm{CGLS}^{\pi} (\cdot , \cdot)$ and the linear functional
$F_\mathrm{CGLS}^{\pi}(\cdot)$ are continuous, existence and
uniqueness of solution of Problem $\mathrm{CGLS}_{\pi}$ follow from
(\ref{eq:cp-coerc}) using the main lemma of \cite{BABUSKA71A}.

We have also considered the following non-symmetric auxiliary problem

\vspace{.25cm}
\paragraph{Problem $\mathrm{CGLS}^{ns}_{\pi}$} {\it Find $\{ \bar{\bm{u}}_h,\bar{p}_h\}
\in \mathcal{U}^l_{h0} \times \mathcal{P}^k_{h0}$  such that
\begin{eqnarray*}
B^{\pi ns}_\mathrm{CGLS}(\{\bar{\bm{u}}_h,\bar{p}_h\},\{\bm{v}_h,q_h\})=
F_\mathrm{CGLS}(\{\bm{v}_h,q_h\}) \quad \forall \; \{{\bm{v}}_h, {q}_h \} \in
\mathcal{U}^l_{h0} \times \mathcal{P}^k_{h0}
\end{eqnarray*}
with
\begin{eqnarray*}
B^{\pi ns}_\mathrm{CGLS}(\{\bar{\bm{u}}_h,\bar{p}_h\},\{\bm{v}_h,q_h\})=
B_\mathrm{CGLS}^1(\{\pi\bm{u}^1_h,p_h\},\{\bm{v}_h,q_h\})&+& \\
B_\mathrm{CGLS}^2(\{\pi\bm{u}^2_h,p_h\},\{\bm{v}_h,q_h\}) &&
\end{eqnarray*}
and $ F_\mathrm{CGLS}(\{\bm{v}_h,q_h\})=
F^1_\mathrm{CGLS}(\{\bm{v}_h,q_h\})+F^2_\mathrm{CGLS}(\{\bm{v}_h,q_h\})$.
}

The asymmetry of this formulation comes from the fact that the
linear transformation $\pi$ is not applied to the weighting
functions as in {Problem $\mathrm{CGLS}_{\pi}$}. We should observe
that this nonsymmetric formulation has the same stability property
of the symmetric one, since, for a given pair
$\{\bar{\bm{u}}_h,\bar{p}_h\}$, we can always choose $\bm{v}_h =
\bar{\bm{u}}_h$ and $ q_h = -\bar{p}_h$ to prove a stability condition
equivalent to (\ref{eq:cp-coerc}). A complete numerical analysis of
these formulations will be presented in a forthcoming paper.
%
%
%
\subsection{Checking the Rates of Convergence }
%
%

\begin{figure}[htb]
\centering
\subfloat[]{\includegraphics[angle=0,width=.49\textwidth,angle=0]{./figs/vyanalitica0-eps-converted-to.pdf}
\label{fig:vy-anis-exata2}}
%
%
\subfloat[]{\includegraphics[angle=0,width=.49\textwidth,angle=0]{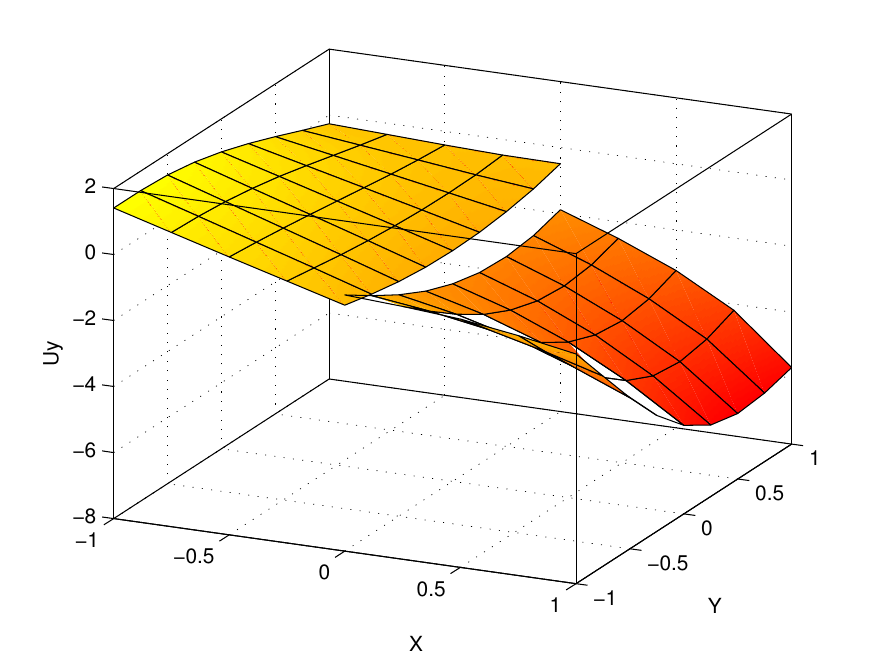}
\label{fig:vy-anis-gls}}
\caption{Component $u_y$;
(a) Exact and (b) approximated by CGLS with $8\times8$ bilinear elements.}
\end{figure}

Now we proceed the convergence study for the anisotropic
heterogeneous problem, considered in Section \ref{sec:numerical1},
by imposing the interface conditions on the stabilized methods
presented before. The exact $y$ component of the velocity field and
its approximation by CGLS with $8 \times 8$ bilinear quadrilaterals
are shown in Figures \ref{fig:vy-anis-exata2} and
\ref{fig:vy-anis-gls}, respectively.
\begin{figure}[htb]
\centering
\subfloat[]{\includegraphics[angle=0,width=.49\textwidth,angle=0]{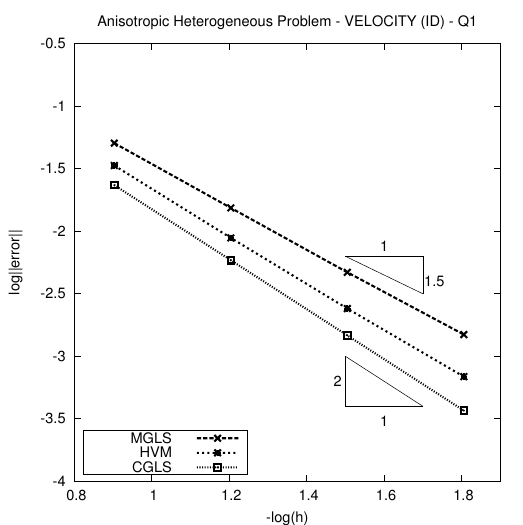}
\label{fig:vel-q1-descontinuo}}
%
%
\subfloat[]{\includegraphics[angle=0,width=.49\textwidth,angle=0]{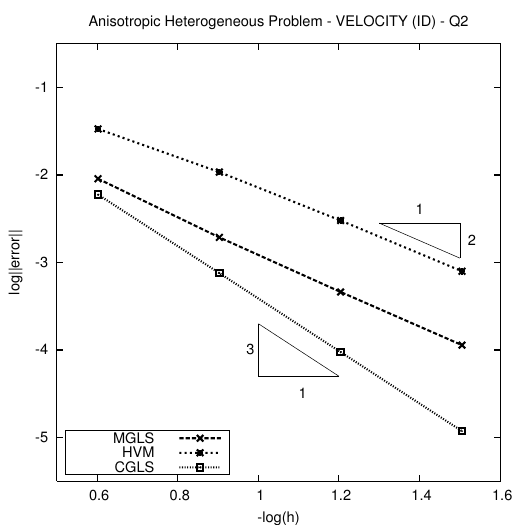}
\label{fig:vel-q2-descontinuo}}
\caption{Convergence study for an anisotropic heterogeneous problem by capturing the discontinuity;
Velocity with (a) bilinear elements and (b) biquadratic elements.}
\end{figure}

Convergence results of velocity approximations for Q1 and Q2
elements with continuous/discontinuous stabilizations are presented
in Figures \ref{fig:vel-q1-descontinuo} and
\ref{fig:vel-q2-descontinuo}, respectively.
The CGLS and HVM methods present convergence rates close to
$O(h^{2.0})$ for Q1, while MGLS leads to $O(h^{1.5})$ convergence.
The results for Q2 elements fits with that expected for homogeneous
problems with regular solutions: $O(h^{2.0})$ for MGLS and HVM, and
to $O(h^{3.0})$ for CGLS.
\begin{figure}[htb]
\centering
\subfloat[]{\includegraphics[angle=0,width=.49\textwidth,angle=0]{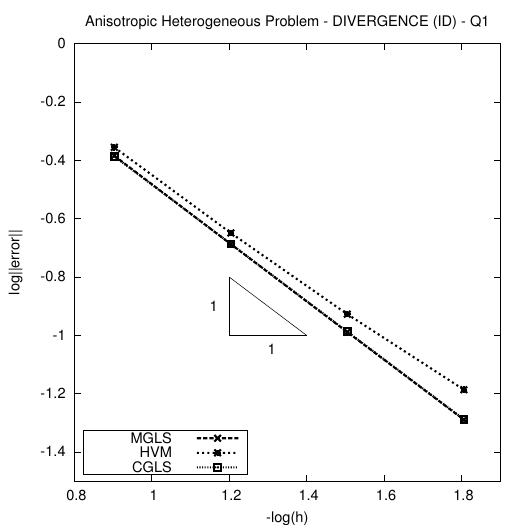}
\label{fig:div-q1-descontinuo}}
%
%
\subfloat[]{\includegraphics[angle=0,width=.49\textwidth,angle=0]{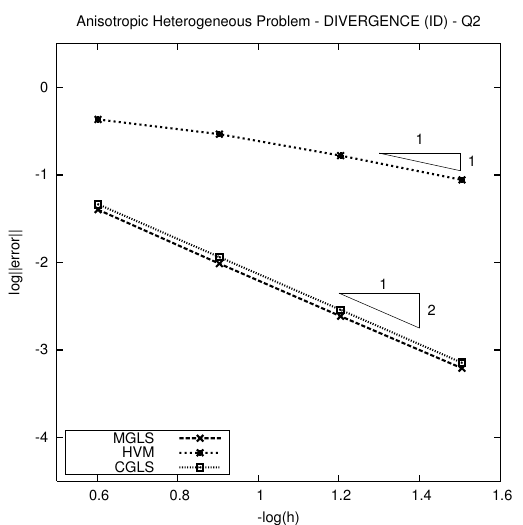}
\label{fig:div-q2-descontinuo}}
\caption{Convergence study for an anisotropic heterogeneous problem
by capturing the discontinuity;
Divergence of velocity with (a) bilinear elements and (b) biquadratic elements.}
\end{figure}

\begin{figure}[htb]
\centering
\subfloat[]{\includegraphics[angle=0,width=.49\textwidth,angle=0]{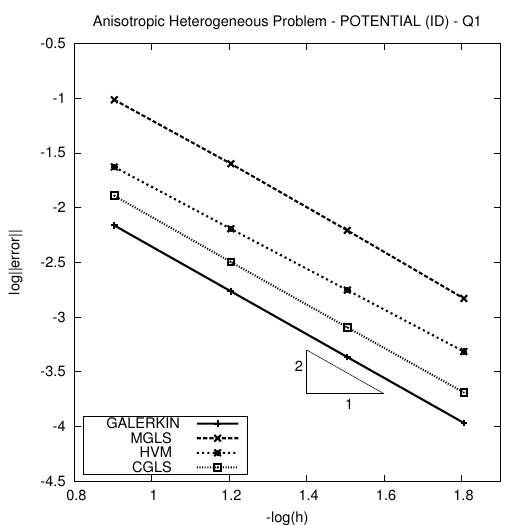}
\label{fig:pot-q1-descontinuo}}
%
%
\subfloat[]{\includegraphics[angle=0,width=.49\textwidth,angle=0]{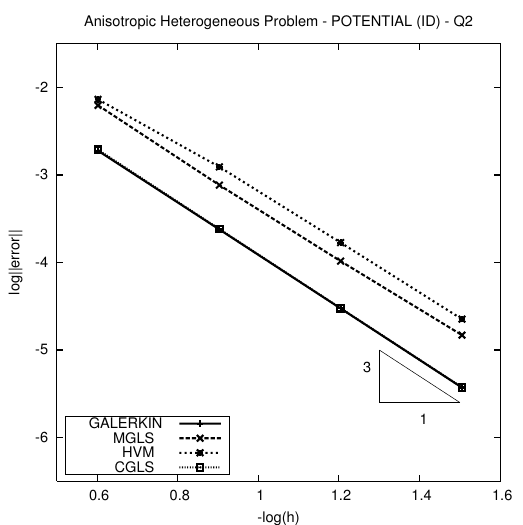}
\label{fig:pot-q2-descontinuo}}
\caption{Convergence study for an anisotropic heterogeneous problem
by capturing the discontinuity;
Potential with (a) bilinear elements and (b) biquadratic elements.}
\end{figure}

The results for the divergence of the velocity, shown in Figures
\ref{fig:div-q1-descontinuo} and \ref{fig:div-q2-descontinuo}, and
for the potential, presented in Figures \ref{fig:pot-q1-descontinuo}
e \ref{fig:pot-q2-descontinuo} indicate the same stability behavior
of the velocity approximation, leading to convergence rates typical
of  those predicted by the numerical analysis for homogeneous
problems with regular solutions.

%
%
\section{Concluding Remarks}
\label{sec:conclu}
%
%

Taking as a model problem Darcy flow in heterogeneous porous media,
we present some comments on finite element approximations of second
order elliptic problems with variable coefficients presenting an
interface of discontinuity.

The starting point is a single field minimization problem associated
with an extended functional proposed by Bevilacqua, Feij\'oo and
Rojas \citep{BEVI74}. We observe that this approach is closely
related to discontinuous Galerkin finite element formulations.

Residual based stabilized mixed methods are reviewed in the context
of Darcy flow in homogeneous porous media and extended to
heterogeneous media with an interface of discontinuity. Stabilized
mixed methods associated with $C^0$ Lagrangian interpolation present
highly stable and accurate approximations for second order elliptic
problems with globally smooth solutions.

As continuous interpolations are not adequate to approximate
discontinuous fields in heterogeneous media, we propose
generalizations of the stabilized mixed methods, based on Galerkin
Least-Squares formulations, for Darcy flow in heterogeneous porous
media with interfaces of discontinuities. By exactly imposing the
continuity/discontinuity constraints on the interface we can recover
stability, accuracy and the optimal rates of convergence of
classical Lagrangian interpolations when applied to globally smooth
problems.

With the proposed formulation, the continuity/discontinuity
interface constraints are  imposed at element level leading to a
global problem similar to the corresponding  $C^0$ finite element
formulation.

Numerical results illustrate the performance of the proposed
methods. Convergence studies for a heterogeneous and anisotropic
porous medium with a smooth interface of discontinuity confirm the
same rates of convergence predicted for homogeneous problem with
smooth solutions.

\section*{Acknowledgements}

The authors thank to  Funda\c c\~ao Carlos Chagas Filho de Amparo \`a
Pesquisa do Estado do Rio de Janeiro (FAPERJ) and to Conselho Nacional
de Desenvolvimento Cient\'{\i}fico e Tecnol\'ogico (CNPq)
for the sponsoring.

\bibliographystyle{elsart-harv}
\bibliography{onfem}

\end{document}